\documentclass[preprint,12pt,3p]{elsarticle}

\usepackage{graphicx}

\usepackage{amssymb}
\usepackage{tikz,tikz-3dplot}
\usetikzlibrary{backgrounds, shapes}   
 
\usetikzlibrary{automata,positioning}

\usepackage{color}

\journal{Journal of Computational Physics}

\begin{document}

\begin{frontmatter}

\title{Recycling Krylov subspaces for CFD applications and a new 
hybrid recycling solver}

\author[label1,label2]{Amit Amritkar\corref{cor1}}
\address[label1]{Department of Mechanical Engineering, Virginia Tech}
\address[label2]{Department of Mathematics, Virginia Tech}

\cortext[cor1]{I am corresponding author}

\ead{amritkar@vt.edu}

\author[label2]{Eric de Sturler}
\ead{sturler@vt.edu}

\author[label2]{Katarzyna {\'S}wirydowicz}
\ead{kswirydo@vt.edu}

\author[label1]{Danesh Tafti}
\ead{dtafti@exchange.vt.edu}

\author[label3]{Kapil Ahuja}
\ead{kahuja@iiti.ac.in}
\address[label3]{Computer Science and Engineering, IIT Indore, India}
\ead{Worked on this project when graduate student at Virginia Tech. Currently Assistant Professor at IIT Indore, India}

\begin{abstract}

We focus on robust and efficient iterative solvers for the pressure Poisson equation in incompressible Navier-Stokes problems. Preconditioned Krylov subspace methods are popular for these problems, with BiCGStab and GMRES(m) most frequently used for nonsymmetric systems. BiCGStab is popular because it has cheap iterations, but it may fail for stiff problems, especially early on as the initial guess is far from the solution. Restarted GMRES is better, more robust, in this phase, but restarting may lead to very slow convergence. Therefore, we evaluate the rGCROT method for these systems. This method 
recycles a selected subspace of the search space (called recycle space) after a restart.
This generally improves the convergence drastically compared with GMRES(m). Recycling subspaces
is also advantageous for subsequent linear systems, if the matrix changes slowly or is constant. However, rGCROT iterations are still expensive in memory and computation time compared with those of BiCGStab. Hence, we propose a new, hybrid approach that combines the cheap iterations of BiCGStab with the robustness of rGCROT. For the first few time steps the algorithm uses rGCROT and builds an effective recycle space, and then it recycles that space in the rBiCGStab solver. 

We evaluate rGCROT on a turbulent channel flow problem, and we evaluate both
rGCROT and the new, hybrid combination of rGCROT and rBiCGStab on a porous medium flow problem. We see substantial performance gains on both problems.

\end{abstract}


\end{frontmatter}


\section{Introduction}
\label{intro}

The current study focuses on iterative solvers for computational fluid dynamics (CFD) applications on structured grids, solving the Navier-Stokes (N-S) and energy equations. 
For very large problems direct methods are impractical, and even for moderate three-dimensional problems they are typically more expensive than iterative methods in work and storage. This is the case even for systems with many right hand sides; see, e.g., \cite{Wang2007}.

The N-S equations are nonlinear partial differential equations (PDEs) that are commonly discretized in space using the Finite Volume Method (FVM), and occasionally the Finite Element Method (FEM). The pressure-velocity coupling frequently requires special consideration. 
This has given rise to various coupled and decoupled solution algorithms,
where the choice of nonlinear solution algorithm 
leads to particular systems of linear equations. 
In terms of computational cost, the solution of the pressure Poisson equation
is the most critical, and incompressible N-S solvers deal with 
this equation in several ways. 
The pressure Poisson equation is elliptic in nature 
but leads to nonsymmetric matrices for non-Cartesian 
meshes or due to certain boundary conditions. Such systems are typically solved by
preconditioned Krylov subspace methods~\cite{templates} with a 
suitable preconditioner. For symmetric positive definite 
systems on regular meshes, multigrid preconditioners or 
multigrid methods are advised, 
but for nonsymmetric systems
the case is less clear. In this paper, we use relatively simple 
multilevel preconditioners that are particularly suited for implementation 
on graphical processing units (GPUs); see below for further details. 
More general multilevel preconditioners or multigrid methods are interesting 
next steps. However, these are not considered here. 

We solve the N-S equations at discrete time steps to capture the transient 
flow characteristics in unsteady flows or to converge to a steady state solution,
leading to a sequence of linear systems of equations. For these systems, originating from a Newton-type iteration, there is little or no change in the coefficient matrix ($A$), and only the right hand side ($b$) changes for subsequent time steps. In this case, recycling a judiciously selected subspace of the search space typically reduces the number of iterations substantially. 
This is referred to as {\em Krylov subspace recycling}~\cite{parks2006recycling}. There are several algorithms that take advantage of recycling for subsequent systems as well as for restarts of generalized minimal residual (GMRES) type methods, such as the generalized conjugate residual method with inner orthogonalization (GCRO)~\cite{de1996nested,de1996inner}, GCRO with deflated restarting (GCRODR)~\cite{parks2006recycling}, GMRES with deflated restarting (GMRES-DR)~\cite{morgan2002gmres} (only for a single system) and GCRO with optimal truncation (GCROT)~\cite{de1999truncation}. For a sequence of nonsymmetric linear systems,  bi-Lanczos-based recycling algorithms, like recycling biconjugate gradient stabilized (rBiCGStab)~\cite{ahuja2014rbicgstab}, have the advantage of a short-term recurrence, and hence low storage requirements and fewer orthogonalizations compared with GMRES variants.

Krylov subspace recycling for sequences of systems has been sporadically applied to structured grid CFD calculations. Several studies in the area of aerodynamic shape optimization show that Krylov subspace recycling with simplified GCROT(m,k)~\cite{hicken2010simplified} improves convergence~\cite{hicken2010parallel,leung2012aerodynamic}. In a study by Carpenter et al.~\cite{carpenter2010general}, using
an enriched GMRES method with subspace recycling for a long sequence of linear systems for steady convection-diffusion problems and flow over a wind turbine, solver parameter tuning yields improved convergence and even eliminates 
stagnation in some cases. For CFD matrices from the Harwell–-Boeing Sparse Matrix Collection, the GCRO-DR algorithm is accelerated by avoiding small skip angles~\cite{Niu2013131}.  
A study by Meng et al.~\cite{meng2014block} shows that block GCROT(m,k) performs better than block GMRES(m) for most of the Poisson and convection-diffusion problems. There are some cases where recycling has not worked particularly well. A study by Mohamed et al.~\cite{mohamedkrylov} observed poor performance of GCRO-DR compared with GMRES(m) for both restarting and sequences of systems for a variety of aerodynamic flows. 

These previous efforts on Krylov recycling in CFD applications all concentrate  
on variants of GMRES. Here, we extend recycling for CFD applications to 
the bi-Lanczos-based method rBiCGStab, and we propose a new hybrid method
as well. The new hybrid method exploits a simplified and more efficient version of the rBiCGStab algorithm derived in \cite{ahuja2014rbicgstab}. This more efficient version made possible by restricting the type of recycling.  
In this paper, we first demonstrate the advantages of Krylov subspace recycling for a turbulent channel flow problem by comparing the results of BiCGStab, GMRES(m) and recycling GCROT (rGCROT). Next, we introduce a novel method for building the outer recycle space for the rBiCGStab algorithm. We use rGCROT to build the recycle space and then switch to rBiCGStab using the same recycle space for subsequent systems. This hybrid method combines the robustness of rGCROT in the initial phase of the simulation, when the starting residual ($r_0$) is large, with the economy of BiCGStab iterations when $r_0$ drops as the flow develops. We test our hybrid method on the pressure Poisson equation for a porous media flow problem with 2.56 million unknowns. We believe this is the first application of preconditioned Bi-Lanczos based recycling algorithms to CFD.

Below, we use rGCROT(m,k) to denote the main parameters, m and k, giving 
the dimensions of the inner Krylov space (as for GMRES(m)) and the outer recycle space, respectively. In the next section, we briefly outline the CFD applications
and the GenIDLEST package, and we discuss the preconditioners and linear solvers. Section~\ref{sec3} describes our CFD test problems. Results from numerical experiments are discussed in section~\ref{sec4} and conclusions in section~\ref{sec5}.

\section{Methodology}
\label{sec2}
In this section, we present a brief overview of the GenIDLEST package and the CFD applications it is used for, along with details of the linear solvers used. We briefly discuss Krylov subspace recycling and we suggest an improved algorithm for rBiCGStab.

\subsection{GenIDLEST}
\label{subsec1}
The GenIDLEST (\textbf{Gen}eralized \textbf{I}ncompressible \textbf{D}irect and \textbf{L}arge-\textbf{E}ddy \textbf{S}imulations of \textbf{T}urbulence) code used in this study solves the incompressible, time-dependent, Navier-Stokes and energy equations in a generalized structured body-fitted multiblock framework. It has been extensively used in propulsion, biological flows, and energy related applications that involve complex multi-physics flows~\cite{taftigenidlest}. GenIDLEST has various turbulence modeling capabilities including Reynolds Average Navier-Stokes (RANS), Detached Eddy Simulation (DES) and Large Eddy Simulation (LES) with subgrid stresses modelled using the dynamic Smagorinsky eddy viscosity model~\cite{tafti1991numerical}. The code uses a finite volume formulation with fractional step algorithm using the semi-implicit Adams-Bashforth/Crank-Nicolson method or a fully-implicit Crank-Nicolson method for a predictor step. The corrector step solves a pressure Poisson equation to satisfy mass continuity~\cite{tafti2009time}. The non-dimensional pressure Poisson equation in generalized coordinates has the following form in GenIDLEST~\cite{tafti2009time},
\begin{equation}
\label{poisson}
\frac{\partial}{\partial \xi_j}\bigg( \sqrt{g}g^{jk} \langle \frac{\partial}{\partial \xi_k} p^{n+1} \rangle \bigg)=\frac{1}{\Delta t} \frac{\partial \langle \sqrt{g} \tilde{U}^j\rangle}{\partial \xi_j} ,
\end{equation}
where  $\sqrt{g}$ is the Jacobian of the spatial transformation, $g^{jk}$ is the contravariant metric tensor, $p^{n+1}$ is the pressure at time level $n+1$ (with density absorbed in the pressure), $\Delta t$ is the time step and $\langle \sqrt{g} \tilde{U}^j\rangle$ is the conserved contravariant flux based on an intermediate velocity $\tilde{u}_i$. Using finite volume discretization, this Poisson equation is transformed into a linear system, $Ax=b$, with pressure as the unknown. The pressure coefficient matrix ($A$) consists of geometric quantities from the left hand side of~(\ref{poisson}). The right hand side ($b$) is evaluated from a local balance of the intermediate volume fluxes at cell faces surrounding a finite volume. The current study focuses on the solution of linear system arising from (\ref{poisson}).

The GenIDLEST code spans over 300,000 lines and more than 600 subroutines, and 
it has modules for the Arbitrary Lagrangian-Eulerian (ALE) Method~\cite{gopalakrishnan2009parallel}, the Discrete Element Method (DEM)~\cite{amritkar2014efficient}, the Immersed Boundary Method (IBM)~\cite{nagendra2014new}, and for Fluid Structure Interaction (FSI)~\cite{gopalakrishnan2010effect}. The computational algorithms 
are optimized to take maximum advantage of state-of-the-art hierarchical 
memory and parallel architectures, with a focus on Message Passing Interface (MPI) and OpenMP-based
codes for central processing units (CPUs)~\cite{Amritkar2012501} and more recently for GPUs~\cite{amritkar2016gpu}.

\subsubsection*{Data Layout}
GenIDLEST uses a non-staggered grid topology; so, all the unknowns are defined and calculated at the computational cell center except the contravariant volume fluxes, which are defined and calculated at cell faces. The coefficient matrix is stored in structured banded sparse storage format (also diagonal storage) \cite[section 11.5]{saad2003iterative}. For every unknown, only the coefficients for the 7 or 19 point stencil, for orthogonal
or non-orthogonal grids respectively, are stored. Although this format (potentially) stores coefficients more than once, it avoids building the stencil more than once and has contiguous memory storage patterns for efficient utilization of CPU cache memory and vectorization. 
If orthogonal meshes are used, and the domain decomposition leads to matching interfaces at mesh block boundaries, the matrices 
from (\ref{poisson}) are symmetric positive definite and conjugate gradient method (CG) is well-suited. In all other cases, 
the system matrices are nonsymmetric and nonsymmetric linear solvers must be used. 
The sparsity pattern of the matrices depends on the boundary conditions.
In GenIDLEST, the matrix from the pressure Poisson equation changes only when there is ALE/IBM grid movement or for IBM related boundary flux correction.

\subsubsection*{Preconditioners in GenIDLEST}
Preconditioners are essential for Krylov subspace methods, as they 
greatly improve the rate of convergence for many problems. In some cases, they are required for convergence. Preconditioners are the most expensive operation in the GenIDLEST framework. Thus a good choice for the preconditioner is imperative.

GenIDLEST uses domain decomposition-based preconditoners that have been optimized for parallel~\cite{wang1999performance} execution on CPUs and, more recently, on GPUs \cite{amritkar2014gpu,amritkar2016gpu}. Optimizations for serial implementations have also been done \cite{wang1998uniprocessor}, but here we focus on preconditioners that run fast in parallel.
In particular for high performance on GPUs, it is important to have a preconditioner tuned for the architecture~\cite{kasia2015}.
For preconditioning, we partition the domain into many overlapping sub-domains, and we use an additive or multiplicative Schwarz method. 

The subdomain equations are solved inexactly using a few 
(local) iterations of Jacobi or SSOR methods. The substructured overlapping subdomain blocks are well-suited to exploit SIMD parallelism and the hierarchical memory of CPUs.
Global coupling between distant subdomains is provided by 
Jacobi iteration on the global problem 
\cite{wang1999performance}. 

Since the forward- and backward substitution in SSOR is less suited to GPUs, we focus on Jacobi iterations for preconditioning in this paper. 
Five sweeps of the Jacobi method with a relaxation parameter (under-relaxation) are applied on the local sub-domains, followed by a sweep of point Jacobi smoothing on the global domain. This preconditioner is referred to as the~\emph{Jacobi preconditioner} in this paper. Additional details of the preconditioners can be found in Wang et al~\cite{wang1999performance}.
The preconditioners are applied from the left of the coefficient matrix, except for BiCGStab, where right preconditioning is used.

\subsection{Linear Solvers}
\label{linsol}
In GenIDLEST, linear solvers are needed for the solution of decoupled momentum equations and the pressure equation. We need to solve three
linear systems in 2D or four in 3D for every time step. CG is used for symmetric systems, and 
BiCGStab or GMRES(m) are used for non-symmetric systems. 
In this paper, we introduce and test in GenIDLEST the recycling Krylov subspace solvers rGCROT~\cite{parks2006recycling,de1999truncation} and rBiCGStab~\cite{ahuja2014rbicgstab}, 
with a new hybrid approach for the latter. 
In the applications considered here, the momentum equations converge rapidly ($\le 2$ iterations) and are always solved with BiCGStab. However, the pressure Poisson equation is more difficult to solve, and we consider various linear solvers described in this section. 
Since the quantity of interest is the pressure correction from (\ref{poisson}), at the start of each time step an initial 
zero pressure field, $x_0=0$, is used for all the solvers. 
Thus, the initial residual is $r_0=b-Ax_0=b$ for every time step.

\subsubsection*{GMRES/GMRES(m)}
\label{gmres}

\begin{figure}[t] \label{fig:GMRES_alg}
{\bf Algorithm 1:} {\it GMRES(m)} \vspace{-3ex}
\begin{tabbing}
......\=....\=....\=....\=....\=....\=.... \kill \\
1.  \> Choose $x_0$, compute $r_0 = b-Ax_0$, $\|r_0\|_2$, and set $i=0$.  \\
2.  \> Choose {\tt tol}, {\tt $m$}, and {\tt max\_itn}.  \\
3.  \> {\bf while} $\|r_i\|_2 > \mathtt{tol} * \|b\|_2$ 
         {\bf and } $i \leq \mathtt{max\_itn}$  \\
4.  \> \> $v_1 = r_i/\|r_i\|_2$  \\
5.  \> \> {\bf for} $j = 1 \ldots m$  \\
6.  \> \> \> $v_{j+1} = Av_j$; \, $i = i+1$  \\
7.  \> \> \> {\bf for} $k = 1 \ldots j$ {\bf do}  \\
8.  \> \> \> \> $h_{k,j} = v_k^T v_{j+1}$; \, $v_{j+1} = v_{j+1} - v_k h_{k,j}$  \\
9.  \> \> \> {\bf end}  \\
10. \> \> \> $h_{j+1,j} = \|v_{j+1}\|_2$; \, $v_{j+1} = h_{j+1}^{-1} v_{j+1}$\\
11. \> \> {\bf end}  \quad \hspace{3.22in}  
              \{\em\, Arnoldi: $AV_m = V_{m+1}\underbar{H}_m$ \}  \\
12  \> \> {\bf Solve} $y = \arg \min_{\tilde{y} \in \mathbb{R}^m} 
            \left\|e_1 \|r_0\|_2 - \underbar{H}_{m} \tilde{y} \right\|_2$ 
            \quad \hspace{1.1in} \{\em\, where $\underbar{H}_m = (h_{k,j})$ \}\\
13. \> \> $x_i = x_{i-m} + V_m y$; \quad \hspace{0.85in}
            \{\em\, where $V_m = [v_1 \, v_2 \, \ldots \, v_m]$, 
            analogous for $V_{m+1}$ \} \\
14. \> \> $r_i = b - Ax_i$; \quad \hspace{1.82in} \{\em\, alternatively 
            $r_i = r_{i-m} + V_{m+1} (\underbar{H}_m y) $ \}  \\
15. \> {\bf end}
\end{tabbing}\vspace{-4ex}
\begin{center}
\parbox[t]{6in}{\small {\bf Note:} For ease of exposition, a few simplifications are made in this presentation.
Our actual implementation handles all these cases properly.
First, the algorithm given above does not check convergence in the inner cycle of $m$ iterations. 
An estimate of the residual is available at each iteration, and the algorithm can be stopped
at an arbitrary iteration. 
Second, the algorithm can stop early if at any point $\|v_{j+1}\|_2 = 0$; this is
not checked in the algorithm above.
}
\end{center}
\end{figure}

The Generalized Minimum Residual Method (GMRES)~\cite{saad1986gmres}, 
given in Algorithm~1, is an 
iterative method for linear systems. At iteration $m$, the method minimizes 
the residual norm, 
$\|b - Ax_m\|_2$, over all vectors $x$ in the Krylov subspace, 
$K_m(A,b) = \mathrm{span}\{b, Ab, \ldots, A^{m-1}b\}$, where 
$A$ and $b$ can be considered as the 
preconditioned matrix and right hand side, respectively.
The Arnoldi recurrence, lines~4--11, generates a sequence of orthogonal vectors that span the Krylov space \cite{arnoldi1951principle}.
These vectors satisfy the (standard) Arnoldi recurrence relation
\begin{equation} \label{eq:arnoldi-rec}
  AV_m = V_{m+1}\underbar{H}_m .
\end{equation}
The solution to the least squares system, line~12, provides
an approximate solution by orthogonal projection.
For a nonsymmetric system, each new Arnoldi vector must be explicitly 
orthogonalized against all previous ones. Therefore storage 
increases linearly with iteration count, until convergence
or a restart (line~15), and work increases 
quadratically with iteration count. 
To limit high costs in storage and computation, {\em restarted} versions of the method, referred to as GMRES(m)~\cite{saad1986gmres}, are used. 
However, restarting often results in a significant increase in the 
total number of iterations. 
The storage cost associated with GMRES(m) in our current implementation 
is $(m+1)N+4N$, in addition to the system 
matrix (here $19N$), 
where $N$ is the number of unknowns and $m$ is the restart frequency.

\subsubsection*{BiCGStab}
\label{bicgstab}

\begin{figure}[t!] \label{fig:BiCGStab_alg}
{\bf Algorithm 2:} {\it BiCGStab} \vspace{-3ex}
\begin{tabbing}
......\=....\=....\=....\=....\=....\=.... \kill \\
1.  \> Choose $x_0$, compute $r_0 = b-Ax_0$, $\|r_0\|_2$, and set $i=0$.
       Choose $\tilde{r}$.  \\
2.  \> Choose {\tt tol} and {\tt max\_itn}  \\
3.  \> {\bf while} $\|r_i\|_2 > \mathtt{tol} * \|b\|_2$ 
       {\bf and} $i \leq \mathtt{max\_itn}$  \\

4.  \> \> $\rho = \tilde{r}^Tr_i$  \\
5.  \> \> {\bf if} $\rho == 0$ {\bf then} breakdown occurred; 
          {\bf exit} gracefully   \\
6.  \> \> {\bf if} $i == 0$ {\bf then} $p = r_i$ {\bf else} 
          $\beta = (\rho / \rho_\mathrm{old})(\alpha/\omega)$;\, 
          $p = r_i + \beta(p - \omega v)$ {\bf end}  \\
7.  \> \> $v = Ap$; \\
8.  \> \> $\alpha = \rho / (\tilde{r}^Tv)$; $s = r_i - \alpha v$;  \\
9.  \> \> {\bf if} $\|s\|_2 \leq \mathtt{tol} * \|b\|_2$ {\bf then} \\
10. \> \> \> $x_{i+1} = x_i + \alpha p $;\, 
             $r_{i+1} = s$;\, 
             {\bf exit} triumphantly (converged)  \\
11. \> \> {\bf end} \\
12. \> \> $t = A s$; \\
13. \> \> $\omega = (t^Ts) / (t^Tt)$  \\
14. \> \> $x_{i+1} = x_i + \alpha p + \omega s$;  \\
15. \> \> $r_{i+1} = s - \omega t$;  \\
16. \> \> $\rho_\mathrm{old} = \rho$; $i = i+1$; \\
17. \> {\bf end}
\end{tabbing} \vspace{-4ex}
\begin{center}\parbox[t]{6in}{
\small{{\bf Note:} Typically, $\tilde{r}$ is either chosen as a 
random vector or equal to $r_0$.}}
\end{center}
\end{figure}

Unlike GMRES, BiCGStab does not need to store the 
entire Krylov subspace, nor does it perform a full orthogonalization 
(Gram-Schmidt process). Instead, BiCGStab relies implicitly
on a (non-optimal) oblique projection to define its iterations~\cite{SimSzy10}, 
which requires only a short recurrence. The BiCGStab method was developed 
by Henk A. van der Vorst~\cite{van1992bi} and quickly gained popularity 
due to its robustness and low computational cost.
The BiCGStab algorithm is given in Algorithm~2.
The short term recurrence in BiCGStab, in iteration $i$,
builds the residual as the product of two polynomials, 
\begin{eqnarray} \label{eq:bicgstab-res}
  r_i = Q_i(A)P_i(A)r_0 .
\end{eqnarray}
As a result, BiCGStab requires two (preconditioned) matrix-vector
products per iteration, with the possible exception of 
the last iteration if convergence is attained at 
line~9 (Algorithm 2).
The storage cost of BiCGStab is $8N$, in addition to the system 
matrix (here $19N$).

\subsection{Recycling Krylov Subspaces}

\begin{figure}[t] \label{fig:GCROT-DR_alg}
{\bf Algorithm 3:} {\it rGCROT(m,k)}/GCRODR(m,k) \vspace{-3ex}
\begin{tabbing}
......\=....\=....\=....\=....\=....\=.... \kill \\
1.  \> Choose $\hat{x}$; compute $\hat{r} = b-A\hat{x}$, 
       $\|b\|_2$, and set $i=0$.  \\
1a. \> Choose/Given $U$, compute $\tilde{C} = AU$ and
       $\tilde{C} = CR$ (thin QR decomposition) \\
1b. \> $\xi = C^T\hat{r}$;\, 
       $r_0 = \hat{r} - C \xi$;\, $\rho_0 = \|r_0\|_2$;\, 
       $x_0 = \hat{x} + U (R^{-1} \xi) $ \\
2.  \> Choose {\tt tol}, {\tt $m$}, and {\tt max\_itn}.  \\
3.  \> {\bf while} $\rho_i > \mathtt{tol} * \|b\|_2$ 
         {\bf and } $i \leq \mathtt{max\_itn}$  \\
4.  \> \> $v_1 = r_i/ \rho_i$  \\
5.  \> \> {\bf for} $j = 1 \ldots m$  \\
6.  \> \> \> $v_{j+1} = Av_j$; \, $i = i+1$  \\
6a. \> \> \> {\bf for} $\ell = 1 \ldots k$ {\bf do} \\
6b. \> \> \> \>  $b_{\ell,j} = c_{\ell}^T v_{j+1}$; \, 
                 $v_{j+1} = v_{j+1} - c_{\ell} b_{\ell,j}$  \\
6c. \> \> \>{\bf end} \quad \hspace{3.2in} 
              \{\em\, $v_{j+1} = v_{j+1} - CC^Tv_{j+1}$ \}  \\
7.  \> \> \> {\bf for} $\ell = 1 \ldots j$ {\bf do}  \\
8.  \> \> \> \> $h_{\ell,j} = v_{\ell}^T v_{j+1}$; \, 
                $v_{j+1} = v_{j+1} - v_{\ell} h_{\ell,j}$  \\
9.  \> \> \> {\bf end}  \\
10. \> \> \> $h_{j+1,j} = \|v_{j+1}\|_2$; \, 
             $v_{j+1} = h_{j+1,j}^{-1} v_{j+1}$  \\
11. \> \> {\bf end} \quad \hspace{1.95in}
          \{\em\, augmented Arnoldi: $AV_m = CB + V_{m+1}\underbar{H}_m$ \} \\
12  \> \> {\bf Solve} $y = \arg \min_{\tilde{y} \in \mathbb{R}^m} 
            \|e_1 \rho_{i-m} - \underbar{H}_{m} \tilde{y} \|_2$  \\
12a.\> \> $z = -By$  \quad \hspace{3.5in} \{\em\, where $B = (b_{\ell,j})$ \}  \\
13. \> \> $x_i = x_{i-m} + V_m y$  
          \\
13a.\> \> $x_i = x_i - U (R^{-1} z)$  \\
14. \> \> $r_i = b - Ax_i$;\, 
          $\rho_i = \|r_i\|_2$ \, 
          \quad \hspace{0.95in} \{\em\, alternatively 
            $r_i = r_{i-m} + V_{m+1} (\underbar{H}_m y) $ \}  \\
14a.\> \> {\bf Update $U$ and $C$ if desired} \quad \hspace{1.75in}
          \{\em\, for details see \cite{de1999truncation,parks2006recycling} \}  \\
15. \> {\bf end}
\end{tabbing} \vspace{-4ex}
\begin{center}
\parbox[t]{6in}{
\small{{\bf Note:} The changes from the GMRES(m) algorithm are given in 
lines~1a--1b, 6a--6c, 12a, 13a, and 14a.}}
\end{center}
\end{figure}

When solving a sequence of slowly changing linear systems,
$A_j x_j = b_j$, the Krylov
subspaces $K_m(A_j,b_j)$ typically contain smaller subspaces 
that are all close to each other and can be approximated by some {\em recurring}
subspace (which need not be a Krylov space).  
The idea of Krylov subspace recycling is (1) to compute this 
recurring subspace efficiently while solving
subsequent linear systems and (2) to improve the 
rate of convergence of the iterative solves  
by iterating orthogonally to this space.

To get an idea of the effectiveness of this approach,
consider restarted GMRES for a single 
system. If GMRES(m) restarts while little progress
has been made, the method restarts with almost the same residual as before. 
Hence GMRES(m) explores a search space very close to 
the search space from the previous $m$ steps. 
Since it has already computed the optimal solution 
from that space, this leads to
very slow convergence or even stagnation. 
By keeping a judiciously selected subspace 
after restarting (or for the next linear system), recycling 
methods compute a
new search space orthogonal to this recycled space, 
which prevents searching near a previous space and often leads to much faster convergence~\cite{parks2006recycling,de1996nested,de1999truncation}.

Next, we discuss how to implement recycling for the GMRES-like methods 
rGCROT and GCRODR and for the short recurrence method BiCGStab.

\subsubsection*{Recycling GMRES: rGCROT and GCRODR}
\label{rgcrot}

The rGCROT(m,k) (or GCRODR(m,k)) algorithm is given in Algorithm~3.
Let $U \in \mathbb{R}^{N \times k}$ and $\mathrm{range}(U)$ 
be a subspace we want to recycle. If $U$ has been updated
or the matrix has changed, we compute the (thin) QR-decomposition 
\[
  CR = AU ,
 \]
where $C \in \mathbb{R}^{N \times k}$ has orthonormal columns and
$R \in \mathbb{R}^{k\times k}$ is upper triangular; see line~1a. 
In line~1b, we orthogonalize the initial residual $\hat{r}$ against $C$,
and update the initial (solution) guess accordingly.
This ensures 
\[
  r_0, v_1 \perp C.
\]

After this initialization, recycling is implemented
by an augmented Arnoldi recurrence, lines~4--11 in 
Algorithm 3, 
the solution of a modified least squares problem,
lines~12--12a,
and a modified update for the approximate 
solution, lines~13--13a.
Otherwise, the method proceeds like GMRES(m);
compare corresponding lines in Algorithms~1 and 3.
Lines 6a -- 6c extend the standard Arnoldi recurrence
given in lines 4--11 of Algorithm 1,
orthogonalizing the new $v_{j+1}$
against $C$ before orthogonalizing against
$v_1, \ldots, v_j$. 
This gives 
\begin{equation} \label{eq:V_perp_C}
  V_{m+1} \perp C
\end{equation} 
and the augmented Arnoldi 
recurrence relation, see (\ref{eq:arnoldi-rec}) for comparison, 
\begin{equation} \label{eq:aug-arnoldi-rec}
  AV_m = CB + V_{m+1}\underbar{H}_m ,
  \qquad \mbox{where } \quad B = C^TAV_m .
\end{equation}
Hence, we approximate  
the residual in a space,  
$\mathrm{range}(V_{m+1}\underbar{H}_m)$, that 
is indeed orthogonal to the 
(previously used) space $\mathrm{range}(C)$.
The optimal update to the solution from the space
$\mathrm{range}(U) + \mathrm{range}(V_m)$ leads to a modified least squares problem.
We have
\begin{equation} 
\label{eq:xm_UV_rm}
  x_m = x_0 + V_m y + UR^{-1}z \quad \Rightarrow \quad
  r_m = r_0 - AV_m y - Cz ,
\end{equation}
then substituting the augmented
Arnoldi recurrence (\ref{eq:aug-arnoldi-rec}) gives
\begin{eqnarray}
  (y,z) & = & \arg \min_{\tilde{y} \in \mathbb{R}^{m}, \tilde{z} \in \mathbb{R}^{k} } 
    \left\| V_{m+1} (e_1 \|r_0\|_2  - \underbar{H}_m\tilde{y}) - 
      C( B\tilde{y} + \tilde{z}) \right\|_2 .
\end{eqnarray}
Using (\ref{eq:V_perp_C}), we can minimize separately for 
the $V_{m+1}$ component (cf. Algorithm 1, line~12)
and the $C$ component. The latter is minimized by taking $z = -By$ (line~12a).

If the recurring subspace approximately persists from one system to the 
next, recycling is very effective. The overhead consists mostly of
the additional orthogonalizations in lines~6a--6c. 
Compared with GMRES(m) recycling variants may even reduce memory
requirements, because $m$ can be much smaller while the 
number of columns in $U$ and $C$ is often modest; see~\cite{de1999truncation}. 
Nevertheless, the memory requirements are typically substantially 
higher than for short recurrence methods like BiCGStab.

GCROT and GCRODR differ only in how they select the space $\mathrm{range}(U)$.
The GCRO-DR algorithm uses approximate invariant subspaces for 
$U$, whereas GCROT measures angles between successive search 
spaces to approximate the recurring subspace directly
\cite{parks2006recycling}. We refer to GCROT with recycling 
for a sequence of systems as rGCROT here.
For our applications of interest, rGCROT gives faster convergence 
and is a more robust solver for the initial time steps.
The additional storage cost of the rGCROT solver over GMRES(m) is 
$2kN$, where $k$ is the (maximum) dimension of the recycle space.

\subsubsection*{rBiCGStab}
\label{rbicgstab}

\begin{figure}[t!]
{\bf Algorithm 4:} {\it rBiCGStab with one recycle space (used for the hybrid algorithm)} \vspace{-3ex}
\begin{tabbing}
.......\=....\=....\=....\=....\=....\=.... \kill \\
1.  \> Choose $\hat{x}$, compute $\hat{r} = b-A\hat{x}$, 
       $\|b\|_2$, and set $i=0$.  \\
1a. \> Choose/Given $U$ and possibly $C$ and $R$ (provided by rGCROT for the
       hybrid algorithm).  \\ 
1b. \> {\bf if} $C$ and $R$ not given {\bf then} compute 
          $\tilde{C} = AU$ and
          $\tilde{C} = CR$ (thin QR decomposition) {\bf end}  \\
1c. \> Compute $\eta_1 = C^T\hat{r}$;\, 
       $r_0 = \hat{r} - C \eta_1$;\, 
       and $\xi = - \eta_1$  \quad 
       \{\em\, accumulate $U$-updates for $x$ in $\xi$ \}  \\
2.  \> Choose {\tt tol}, {\tt max\_itn}, and $\tilde{r}$.  \\
3.  \> {\bf while} $\|r_i\|_2 > \mathtt{tol} * \|b\|_2$\, 
         {\bf and } $i \leq \mathtt{max\_itn}$  \\

4.  \> \> $\rho = \tilde{r}^Tr_i$  \\
5.  \> \> {\bf if} $\rho == 0$ {\bf then} breakdown occurred; 
          {\bf exit} gracefully   \\
6.  \> \> {\bf if} $i == 0$ {\bf then} $p = r_i$ {\bf else} 
          $\beta = (\rho / \rho_\mathrm{old})(\alpha/\omega)$;\, 
          $p = r_i + \beta(p - \omega v)$ {\bf end}  \\
7.  \> \> $v = Ap$;  \\
7a. \> \> $\eta_1 = C^Tv$;\, $v = v - C\eta_1$  \\
8.  \> \> $\alpha = \rho / (\tilde{r}^Tv)$;\, $s = r_i - \alpha v$  \\
9.  \> \> {\bf if} $\|s\|_2 \leq \mathtt{tol} * \|b\|_2$ {\bf then}  \\
10. \> \> \> $x_{i+1} = x_i + \alpha p$;\, 
             $r_{i+1} = s$  \\
10a.\> \> \> $\xi = \xi + \alpha \eta_1$;\, 
             {\bf exit} triumphantly (converged)\\
11. \> \> {\bf end}  \\
12. \> \> $t = As$  \\
12a.\> \> $\eta_2  = C^Tt$;\, $t = t - C \eta_2$  \\
13. \> \> $\omega = (t^Ts) / (t^Tt)$  \\
13a.\> \> $\xi = \xi + \alpha \eta_1 + \omega \eta_2$ \quad
          \{\em\, accumulate $U_k$ updates for $x$ \} \\
14. \> \> $x_{i+1} = x_i + \alpha p + \omega s$  \\
15. \> \> $r_{i+1} = s - \omega t$  \\
16. \> \> $\rho_\mathrm{old} = \rho$; $i = i+1$  \\
17. \> {\bf end}  \\
17a.\> $x_i = x_i - U(R^{-1}\xi)$ \quad
       \{\em\, add accumulated updates to solution \}
\end{tabbing} \vspace{-4ex}
\begin{center}
\parbox[t]{6in}{
\small{ {\bf Note:} 
As the solution is needed only at the very end, all updates
with the $U$ vectors (updates of the type $x = x - UR^{-1}z$) are
postponed, reducing $\# \mathrm{its} * 4Nk$ work to $2Nk$ work. The 
length $k$ vector $\xi$ accumulates these postponed updates.
The changes from the standard BiCGStab algorithm are given
in lines~1a--1c, 7a, 10a, 12a, 13a, and 17a.}}
\end{center}
\end{figure}

As we saw above, BiCGStab builds polynomials of the type
(\ref{eq:bicgstab-res}).  
In general, for recycling in BiCGStab
we would recycle two subspaces, each augmenting
one of the two Krylov subspaces in the 
underlying bi-Lanczos recurrence~\cite{ahuja2014rbicgstab}. 
However, as shown in~\cite[section 5, example 1]{ahuja2014rbicgstab}, 
this is not always necessary. If only one recycle space
is used, $\mathrm{range}(U)$, the algorithm 
simplifies substantially. 
To see how such a recycling BiCGStab can be implemented,
we rearrange the augmented Arnoldi recurrence
relation (\ref{eq:aug-arnoldi-rec}) for $A$
as a standard Arnoldi recurrence relation
for $A_C = (I-CC^T)A$. Moving the term $CB$ to the left
and using $B = C^TAV_m$, we get
\begin{equation} \label{eq:aug-arnoldi-equiv}
  AV_m = CB + V_{m+1}\underbar{H}_m \quad \Leftrightarrow \quad
  (I-CC^T)AV_m = V_{m+1}\underbar{H}_m. 
\end{equation}
Hence, as for rGCROT, we can generate the 
Krylov subspace $K_m(A_C,r_0)$ 
with the operator $A_C$ and 
with $r_0$ orthogonalized against $C$; see
lines~1a, 7a, and 12a in Algorithm~4. This leads to 
polynomials of the type
\begin{eqnarray} \label{eq:rbicgstab-res}
  r_{i}^{\mathrm{rBiCGStab}} & = & Q_i( A_C ) P_i( A_C )r_0 ,
\end{eqnarray}
with starting residual $r_0 = (I-CC^T)\hat{r}$. 
So, we run BiCGStab 
with a modified 
initial residual, $A_C$ replacing $A$, and 
special updates for the solution as for 
rGCROT. The special updates for the solution can 
be derived from the fundamental relation between solution
updates and residual updates in Krylov methods,
\begin{eqnarray}
  x_{i+1} = x_i + z & \Rightarrow & r_{i+1} = r_i - Az , 
\end{eqnarray}
which implies
\begin{eqnarray} \label{eq:sol-update-from-res}
  z & = & A^{-1}(r_i - r_{i+1}) .
\end{eqnarray}
Note that (\ref{eq:sol-update-from-res}) still holds
when the residual is computed as in (\ref{eq:rbicgstab-res}).
Hence, the residual update $s = r_i - \alpha A_C p$
(lines~7--8) corresponds to the solution update
\begin{equation}
  x_{i+1} = x_i + A^{-1}(\alpha A_c p) = 
  x_i + \alpha p - U (\alpha R^{-1} \eta_1) 
  \quad \mbox{ (see lines~10--10a) } ,
\end{equation}
where the $U$-component of the 
update, $U(\alpha R^{-1} \eta_1)$, is postponed till 
the end; see line~17a. 
The residual update in line~15 leads to an analogous
update for the approximate solution.

This implementation also avoids some expensive updates
if possible.
The algorithm for the general case of rBiCGStab and its 
derivation can be found in~\cite{ahuja2014rbicgstab,ahuja2011recycling}.
The additional storage for rBiCGStab compared with BiCGStab is $2kN$, 
where again $k$ is the (maximum) dimension of the recycle space.

\subsection{Hybrid approach}
\label{hybrid}

\begin{table}[tbp]
\centering
\caption{Advantages and disadvantages of various linear solvers considered.}
\begin{tabular}{l|p{6.5cm}|p{6.5cm}}
\hline
 Solver & Advantages & Disadvantages \\
 \hline
 BiCGStab 
   & \parbox[t]{0.5em}{$\cdot$}
     \parbox[t]{2.3in}{Cheap iterations (in time) with low         storage requirements} 
   & \parbox[t]{0.5em}{$\cdot$}
     \parbox[t]{2.3in}{Irregular and/or slow convergence for        stiff problems, possible accuracy problems} \\ 
   & \parbox[t]{0.5em}{$\cdot$} 
     \parbox[t]{2.3in}{Often fast in run time}   & \\
 \hline
 rBiCGStab 
   & \parbox[t]{0.5em}{$\cdot$}
     \parbox[t]{2.3in}{Relatively cheap iterations (in time)} 
   & \parbox[t]{0.5em}{$\cdot$}
     \parbox[t]{2.3in}{Improved but irregular convergence} \\
   & & \parbox[t]{0.5em}{$\cdot$}
     \parbox[t]{2.3in}{High storage requirement} \\
 \hline
 GMRES(m) 
   & \parbox[t]{0.5em}{$\cdot$}
     \parbox[t]{2.3in}{Min. residual over $m$ iterations}
   & \parbox[t]{0.5em}{$\cdot$}
     \parbox[t]{2.3in}{Large storage requirement} \\
   & \parbox[t]{0.5em}{$\cdot$}
     \parbox[t]{2.3in}{Monotonic residual norm decrease} 
   & \parbox[t]{0.5em}{$\cdot$}
     \parbox[t]{2.3in}{Expensive iterations due to
             orthogonalization} \\
 \hline
 rGCROT 
   & \parbox[t]{0.5em}{$\cdot$}
     \parbox[t]{2.3in}{Fast convergence (in iterations)            because of recycling} 
   & \parbox[t]{0.5em}{$\cdot$}
     \parbox[t]{2.3in}{Higher storage requirement} \\
   & \parbox[t]{0.5em}{$\cdot$}
     \parbox[t]{2.3in}{Monotonic residual decrease} 
   & \parbox[t]{0.5em}{$\cdot$}
     \parbox[t]{2.3in}{Expensive iterations due to 
       additional orthogonalizations} \\
   \hline
  \end{tabular}
\label{tab:method1}
\end{table}

Although the rGCROT method is more robust and often converges 
in fewer iterations than BiCGStab,
its iterations are relatively expensive. 
As a result, for later time
steps, the run time for rGCROT might be higher than for (r)BiCGSTab.  
In our application, the additional robustness and fast convergence 
of rGCROT are important mainly in the first few time steps.
For later time steps, we can use BiCGStab or, 
for faster convergence and improved robustness, rBiCGStab. 
Therefore, we propose a hybrid strategy that starts with rGCROT and 
switches to rBiCGStab, recycling the subspace computed by rGCROT for previous systems. 
The recycle spaces constructed in rGCROT, $\mathrm{range}(U)$ and 
$\mathrm{range}(C)$, remain 
constant during the rBiCGStab calculations. Thus a large recycle space 
can be used efficiently in this approach. An alternative approach of 
intermittent solves with rGCROT to update the recycle space can be used 
when the system matrix changes for certain time steps.
Moreover, if convergence is poor at an intermediate phase, we can switch back to rGCROT. 
This was not necessary in the applications discussed in this paper. 
We list the advantages and disadvantages of the various solvers 
in Table~\ref{tab:method1}. The hybrid approach is designed to 
combine the advantages of the robust rGCROT algorithm and the 
cheaper rBiCGStab iterations.

\section{CFD problem description}
\label{sec3}
In this study, we consider two problems where the pressure coefficient matrix remains constant for all subsequent time steps,
\begin{enumerate}
\item Turbulent channel flow,
\item Flow through porous media.
\end{enumerate}

\subsection{Turbulent channel flow}
\label{turb}
\begin{figure}
    \centering
    \includegraphics[width = 10 cm]{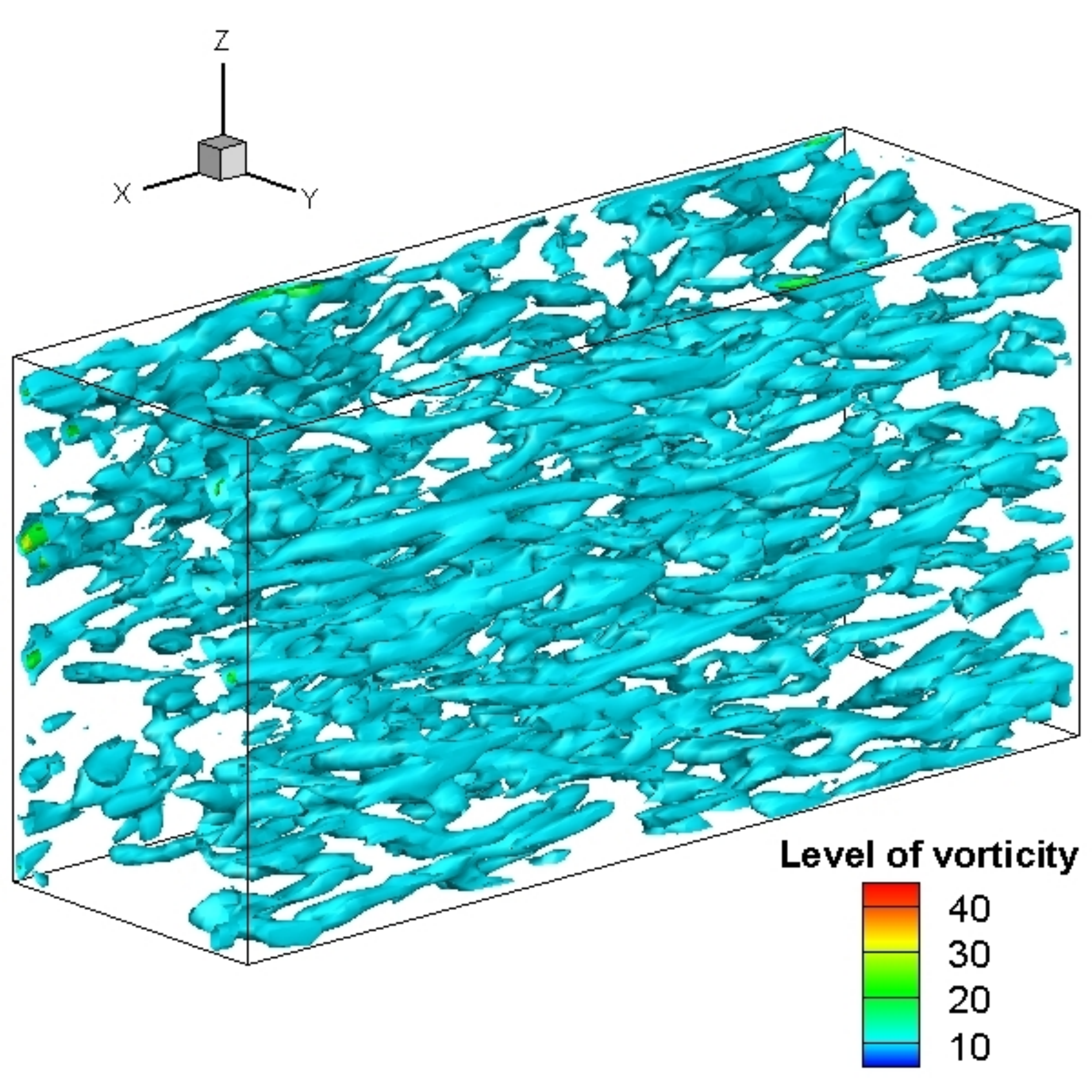}
    \caption{Turbulent channel flow showing vorticity}
    \label{fig:turb}
\end{figure}    

Turbulent channel flow has been used extensively as a canonical problem that embodies most of the physical complexities of wall bounded turbulent flow. For this reason, turbulent channel flow, depicted in Figure~\ref{fig:turb}, is a good problem to evaluate the linear solvers. In this study, we use a 
three-dimensional, orthogonal computational grid,
leading to a seven point stencil, with a single mesh block
of $64 \times 64 \times 64$ computational cells. 
A higher grid density near the wall is used to resolve the boundary layer, and periodic boundary conditions are used in the flow direction to mimic an infinitely long turbulent channel. 
We specify a fixed mean pressure gradient in the flow  direction 
to balance the wall friction. Prescribed perturbations in the flow, based on the mean turbulent channel flow profile, are given as the initial 
conditions to trigger the onset of turbulence. 
The flow is allowed to evolve in time until the solution reaches a stationary state. A skin friction Reynolds number ($Re_\tau$) of 180, based on channel half width and wall friction velocity, is used to compare with the study of Moser et al.~\cite{moser1999direct}.  First, the simulation is run using the rGCROT algorithm until the flow turbulence has become statistically stationary in time (stationary flow). For performance comparison, 30 time steps of $5\times10^{-5}$ second are run at the start and end (after 2.5 seconds) of the calculations with various solvers. An absolute tolerance of $10^{-6}$ is used as the convergence criterion on the $L^2$ norm of the residual. The Jacobi preconditioner with 5 inner iterations is used for all the solvers.

\subsection{Flow through porous media}
\label{porous}
\begin{figure}
    \centering
    \includegraphics[width = 12 cm]{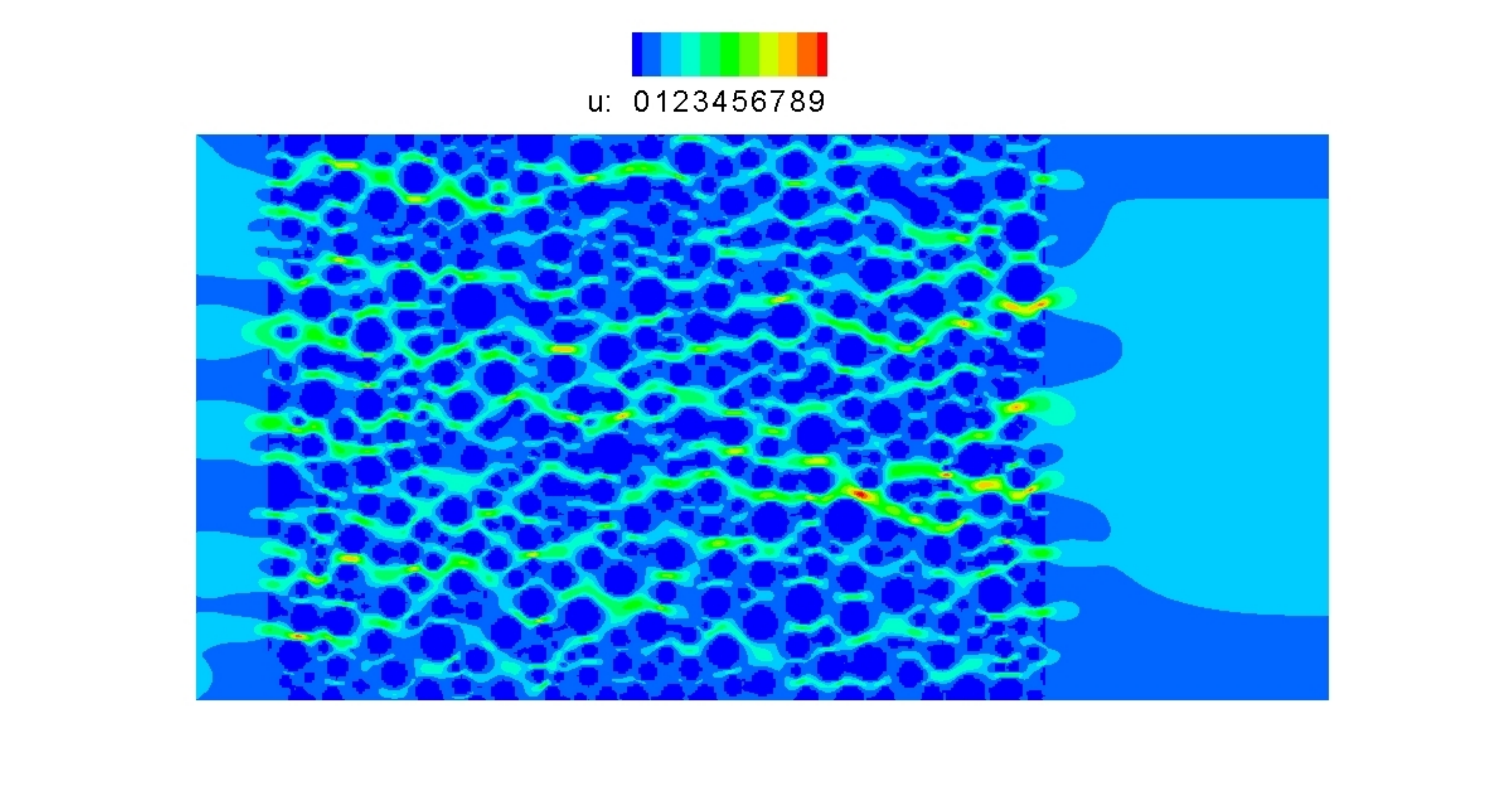}
    \caption{Flow through porous media with $u$ velocity contours}
    \label{fig:porous}
\end{figure}    

Porous media flow has a wide variety of applications. In the current study, 
we model the porous medium using IBM to resolve the porous structure, which is constructed using the stochastic reconstruction procedure, as shown in Figure~\ref{fig:porous}. The IBM scheme uses indirect forcing of sharp interfaces on the background mesh~\cite{nagendra2014new}. 
We use a 2D background mesh consisting of 2.56 million orthogonal 
cells and 16 mesh blocks ($100 \times 800 \times 2 \times 16$) 
to simulate a bulk flow Reynolds number ($Re_b$) of $10^{-4}$. 
In the current IBM framework, two computational cells are needed in the $Z$ direction for 2D calculations. The porous medium has wall boundaries 
on the top and the bottom, with inlet and outlet boundaries on the left and right, respectively. The pressure coefficient matrix remains constant for every time step, as the IBM related flux corrections at the immersed boundaries are not applied. For the performance comparison, 
ten time steps of $10^{-8}$ second are run at the start of the calculations with various solvers. 
A relative tolerance of $10^{-10}$ on the $L^2$ norm of the residual is used as the convergence criterion, as the initial residual norm is large and no solver would converge to an absolute tolerance of $10^{-6}$, since this is too close to the relative machine precision. The preconditioner uses 5 inner iterations of either Jacobi iterative smoothing or SSOR smoothing. Although the SSOR algorithm is more appropriate for symmetric matrices, it can be applied to nonsymmetric coefficient matrices as a preconditioner, 
as this only requires an approximate solution.

\section{Results}
\label{sec4}
The turbulent channel flow simulation is performed serially on a local system (Dual Intel\textregistered Xeon\textregistered CPU X5650 @ 2.67GHz \& 48GB memory), whereas the flow through porous media study is performed using 16 CPU cores (Dual Intel\textregistered Xeon\textregistered CPU E5-2670 @ 2.60GHz \& 64GB memory) with MPI parallelism on the BlueRidge HPC system at Virginia Tech.

In the experiments discussed in this section, we check the convergence for BiCGStab and rBiCGStab only at the end of one complete execution of the while-loop, that is, going once through lines~3--17, in Algorithms~2 and 4. Therefore, one iteration of (r)BiCGStab includes two (preconditioned) matrix-vector products. Similarly, for GMRES(m) and rGCROT(m,k) we check convergence only at the end of one complete execution of the while-loop, that is, going once through lines~3--15, in Algorithms~1 and 3. 
In the context of restarted GMRES(m) and rGCROT(m,k), this is often called a cycle.
A cycle involves $m$ steps of the (augmented) Arnoldi recurrence, and hence 
involves $m$ matrix-vector products. Notice that the number of matrix-vector 
products determines the dimension of the Krylov space from which a solution is
computed. Therefore, in comparing how fast methods converge, we compare the number of matrix-vector products until a specified tolerance is reached.
Since the number of flops (and hence run time) per GMRES(m) or rGCROT(m,k) 
iteration, typically defined as one iteration of the (augmented) Arnoldi recurrence, 
is not constant, we provide some of the timing/performance information below per cycle
($m$ matrix-vector products).
Similarly, we provide some of the timing/performance information for (r)BiCGStab per iteration ($2$ matrix-vector products).

In our comparisons of the various solvers, we have to make a choice of either running each solver independently in GenIDLEST for the specific problem, or choosing one solver as the master and make the solutions computed by that solver determine the linear systems solved by all the other solvers. The first approach has the advantage that the performance numbers reported are the actual results obtained for running GenIDLEST with that solver. It has the disadvantage that the solvers being compared all solve slightly different problems, potentially impacting the comparison. In our case, only the right hand sides are slightly different. The second approach avoids this problem, but will compute performance characteristics that would not actually be observed using these solvers independently for GenIDLEST.

Since we compare the solvers over a fairly large number of time steps 
(linear systems), we prefer the first approach. The effect of a few potentially unfortunate right hand sides for, say, one solver will mostly be averaged out. Moreover, we are careful to ensure that the right hand sides  for each time step don't vary much among the solvers. The starting point for all the solvers is the same, and we use a relatively high convergence criterion to ensure that the solutions of linear systems (which determine subsequent right hand sides), are very close. So, although the right hand sides are not exactly the same across all solvers, we don't expect this to have a large impact on the average and total matrix-vector product counts and on the run time required for convergence. As will be clear from the reported results, for most solvers the number of matrix-vector products for convergence does not fluctuate much from time step to time step; see figures 4 and 5. This suggests the convergence rates are not particularly sensitive to the variations in typical right hand sides.

\subsection{Turbulent channel flow}
\label{res:turb}
There are 6 parameters associated with the rGCROT algorithm \cite{parks2006recycling}: the maximum number of steps of the augmented Arnoldi recurrence (lines~5--11 in Algorithm~3), $m$ 
\footnote{The parameter $m$ also gives the maximum dimension of the inner search space, 
$\mathrm{span}(v_1,\; \ldots\; ,v_m)$, and it is similar to the restart frequency of GMRES.}; the (max) dimension of the recycle space, or outer search space, $k$; the number of outer vectors after truncation of the outer space, here $k-10$; the number of inner vectors to select outer vectors from, here $m/2$; the number of inner vectors selected to extend the the outer space, here $1$; and, the number of latest inner vectors kept to extend the outer space, here $0$.

In this study, we focus on the two most influential parameters, $m$ and $k$. The other parameter values are given as constants or in terms of $m$ and $k$ (as stated above). To analyze optimal performance, $k$ is varied from $20$ to $170$, and $m$ is varied from $20$ to $100$, both with increments of $10$, at the start of the simulations and also at the end of simulations. At the start of the simulations, rGCROT(30,130) is optimal with respect to solution time, whereas rGCROT(20,130) is optimal after the flow becomes stationary. Additional details regarding the rGCROT parameters and selection guidelines can be found in~\cite{parks2006recycling}.

\begin{figure}
    \centering
    \includegraphics[width = 12 cm]{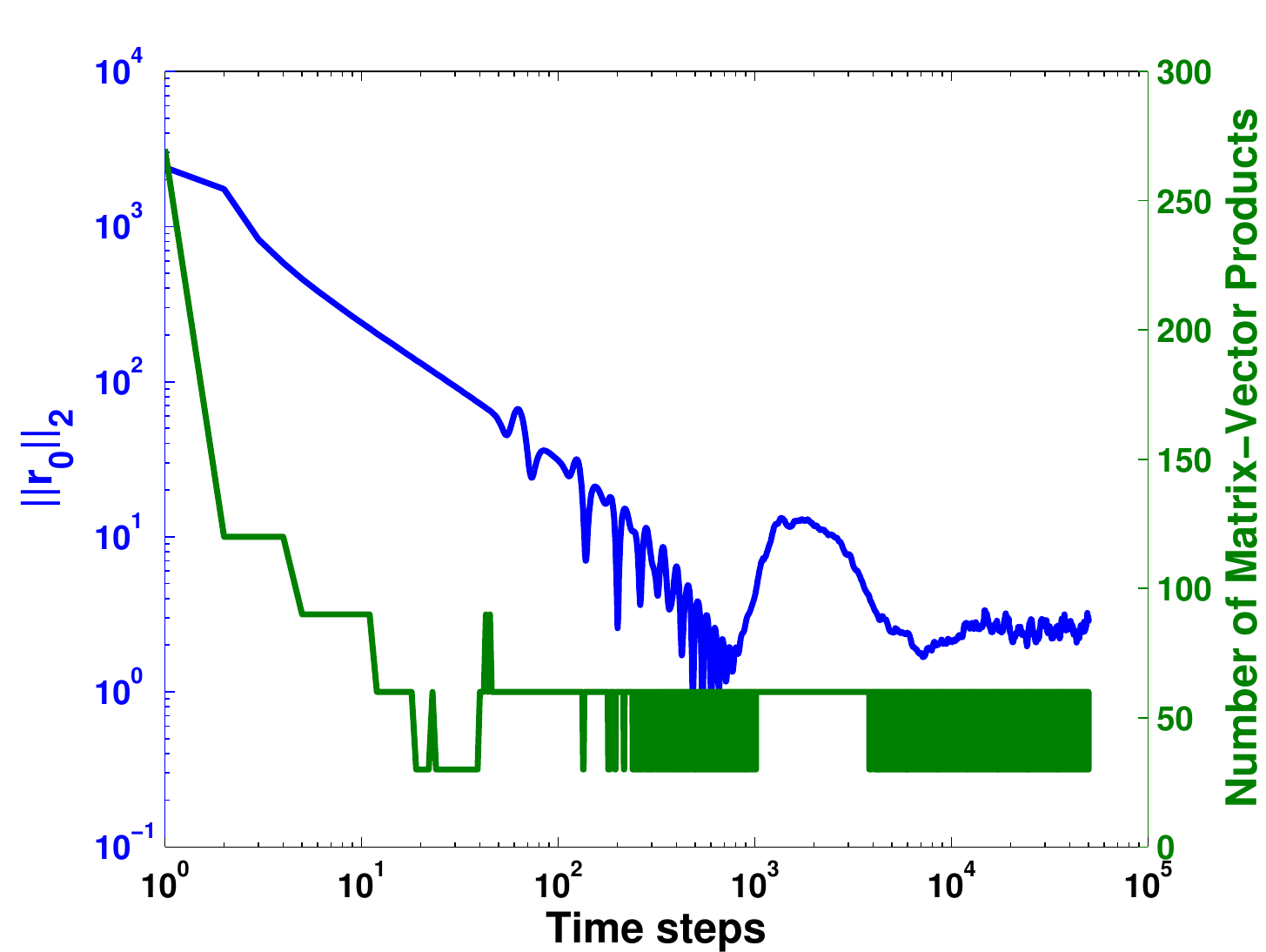}
    \caption{Number of matrix-vector products for rGCROT(30,130)
    to converge for each time step (right y-axis) and 
      the initial $\|r_0\|_2$ for each time step
      (left y-axis).}
    \label{fig:turb1}
\end{figure}

Figure~\ref{fig:turb1} gives the number of matrix-vector products needed to converge for a long sequence of time steps with the rGCROT(30,130) solver. It also shows the $L^2$ norm of the initial residual ($\|r_0\|_2=\|b\|_2$) for every time step. The variation in $\|r_0\|_2$ is directly associated with the flow physics. At the first time step, the initial guess for the velocity with prescribed perturbations in the flow is farthest from the solution, so the $\|r_0\|_2$ value is maximum. As the flow starts to develop, $\|r_0\|_2$ goes down and drops by 3 orders of magnitude. After about 1000 time steps, $\|r_0\|_2$ goes back up by an order of magnitude, because rapid changes in the flow rate cause relatively large changes in subsequent right hand sides ($b$). Once the flow becomes stationary, the flow rate and consequently $\|r_0\|_2$ settle down to quasi-steady values. The number of matrix-vector products needed to converge 
for a given time step is influenced by the variations in $\|r_0\|_2$, since an absolute convergence tolerance of $10^{-6}$ is used. For the rGCROT algorithm, the dimension and quality of the recycle space also influences the number of matrix-vector products needed to converge. Initially, the dimension of the recycle space is small and the recurring subspace has not been constructed yet, and thus the number of matrix-vector
products to converge is higher. As a recurring subspace is discovered and captured by the recycle space, the number of matrix-vector products drops to $30$ or one cycle. 
When the flow develops, the right hand side changes considerably with each 
time step due to rapid changes in flow rate. This warrants intermittent updates to the recycle space, and the number of cycles (each 30 matrix-vector products) to converge fluctuates between one and two. As the $\|r_0\|_2$ values goes back up, the number of cycles increases to two, since the relative reduction required in the residual norm also increases. 
After this phase of fluctuations, the flow settles and the number of
cycles is again one or two.

Tables~\ref{tab:turb1} and~\ref{tab:turb2} give the time to solution and 
the average number of matrix-vector products to converge. 
From the results it is clear that the 
parameter-tuned rGCROT solver takes the least number of
matrix-vector products per time step and has the lowest time to solution. 
Notice that the time for a BiCGStab iteration, 
which takes 2 matrix-vector products and a few orthogonalizations, 
is very small compared with a rGCROT(30,130) cycle, which 
takes 30 matvecs and many orthogonalizations. However, 
it takes many more BiCGStab iterations to converge than 
rGCROT(30,130) cycles. 
At the start of the calculations, GMRES(m) did not converge in 1000 
matrix-vector products if the restart frequency, $m$, was chosen 
less than 50. 
This prolonged stagnation in convergence is a typical 
issue with GMRES(m) for small $m$. 

All the solvers converge faster when the flow becomes stationary, 
because $\|r_0\|_2$ is smaller.

\begin{figure}
    \centering
    \includegraphics[width = 12 cm]{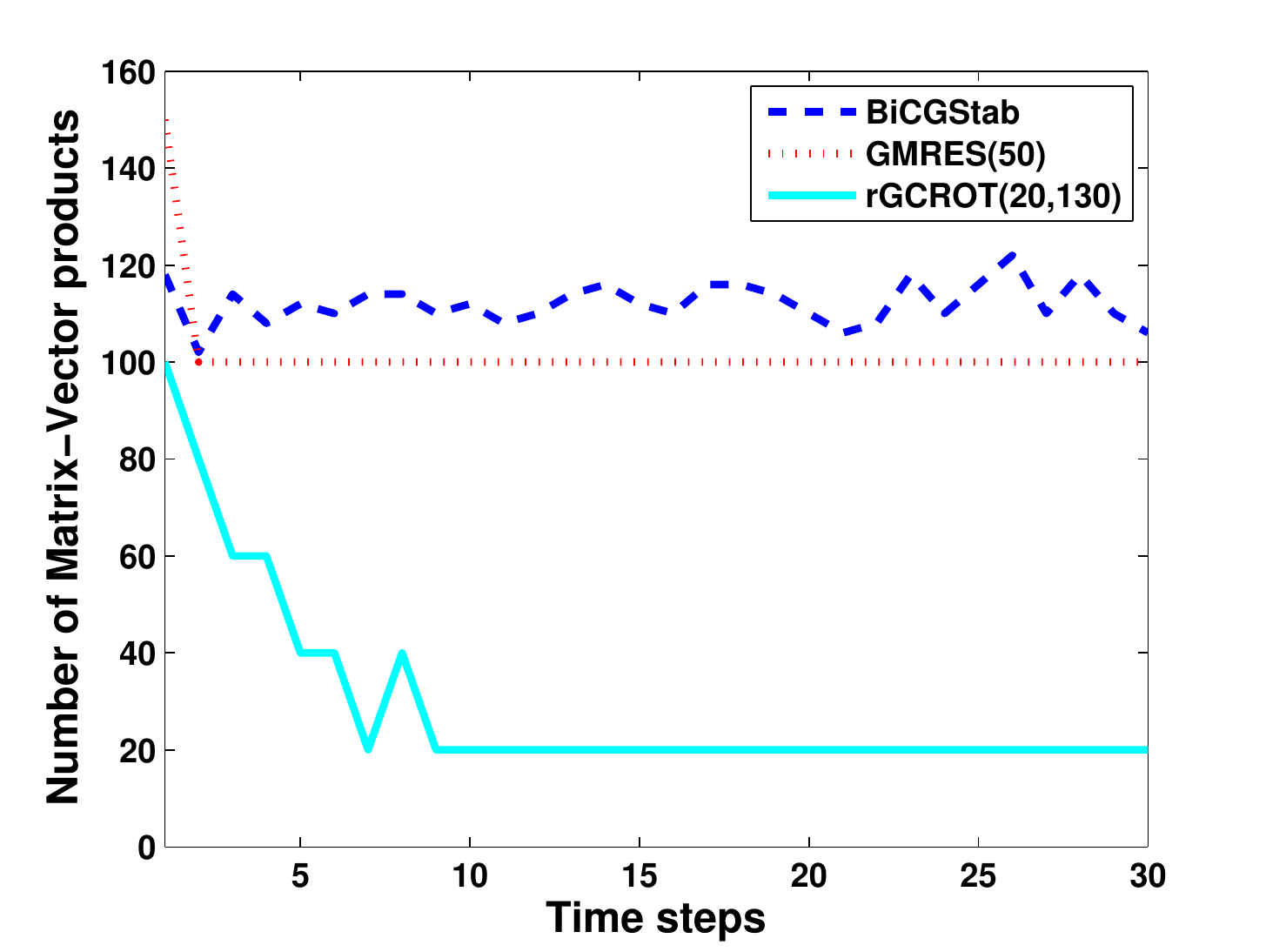}
    \caption{Number of matrix-vector products for each solver 
    to converge to an absolute tolerance of $10^{-6}$ 
    for each time step
    when the flow is stationary.}
    \label{fig:turb2}
\end{figure}

The number of matrix-vector products to converge, after the flow becomes stationary, 
is shown in Figure~\ref{fig:turb2} for a sequence of 30 time steps. The number of
matrix-vector products for rGCROT(20,130) remains constant once the recycle space has captured the dominant subspace. The method convergences in one cycle (20 matvecs). 
The number of matrix-vector products to converge for GMRES(m) 
also remains constant (100 matvecs or 2 cycles), but the time to solution 
is much larger as we have to build the search space from 
scratch after every restart. 
The number of matrix-vector products for BiCGStab fluctuates only modestly.

The storage costs of the various solvers for the turbulent 
channel flow problem are listed in Table~\ref{tab:store1}. 
We only list storage costs beyond the standard costs in every solver, for the system matrix, the right hand side, and the approximate solution, and that are at least O(N), where N is number of unknowns. Required auxiliary and temporary vectors used in the solvers are counted towards the storage overhead. As mentioned in section~\ref{linsol}, the BiCGStab solver has the lowest storage requirements, and the rGCROT solver often has  substantially higher storage requirements. For the same amount of storage, it has been shown that the rGCROT algorithm 
tends to outperform GMRES(m)~\cite{parks2006recycling}. When storage is cheap, 
the rGCROT solver typically shows superior performance compared with other solvers.

In summary, rGCROT is a robust solver that is capable of converging in very small numbers of cycles (or matrix-vector products) by exploiting subspace recycling. However, an rGCROT cycle may be relatively expensive in time and in storage. On the other hand, BiCGStab has cheap iterations but may take many more iterations (or matrix-vector products) to converge. This suggests a hybrid approach that combines the robustness of rGCROT with the cheap iterations of BiCGStab by using the rGCROT recycle space in rBiCGStab as discussed below. We also tested this hybrid approach for the turbulent channel
flow problem, but the method was not faster (in time) than rGCROT.
However, for the porous media application in the next section the hybrid method
does produce shorter runtimes.

\begin{table}[tbp]
\centering
\caption{Solver performance at the start of the calculations}
\begin{tabular}{ p{6cm} | c | c | c }
    \hline  
    First 30 time steps & BiCGStab & GMRES(50) & rGCROT(30,130) \\
    \hline                       
    Total time to solution (s) & 365 & 12680 & 323 \\
    Average number of matrix-vector & 184 & 4350 & 69 \\
    \hspace{1em} products per time step & & & \\
    Average number of iterations (its) & 92 (its) & 87 (c) & 2.3 (c) \\
    \hspace{1em} or cycles (c) per time step & & & \\ 
    Solver time per iteration (s/it)  
       & 0.132 (s/it) & 4.84 (s/c) & 4.68 (s/c) \\
    \hspace{1em}  or per cycle (s/c) & & & \\
    \hline  
\end{tabular}
\label{tab:turb1}
\end{table}

\begin{table}[tbp]
\centering
\caption{Solver performance when the flow is stationary}
\begin{tabular}{ p{6cm} | c | c | c }
    \hline  
    30 time steps at end & BiCGStab & GMRES(50) & rGCROT(20,130) \\
    \hline                       
    Total time to solution (s) &  242.7 & 336.7 & 164.8 \\
    Average number of matrix-vector & 112 & 100 & 30 \\
    \hspace{1em} products per time step & & & \\
    Average number of iterations (its)  & 56 & 2 & 1.5 \\
    \hspace{1em} or cycles (c) per time step & & & \\
    Solver time per iteration (s/it) & 0.144 & 5.52 & 3.66 \\
    \hspace{1em} or per cycle (s/c) & & & \\
    \hline
\end{tabular}
\label{tab:turb2}
\end{table}

\begin{table}[tbp]
\centering
\caption{Solver storage cost comparison for the turbulent channel case.}
\begin{tabular}{l|c}
    \hline
    Solver & Additional storage cost \\
    \hline  
    BiCGStab & $8N$ \\
    GMRES(50) & $55N$ \\
    rGCROT(30,130) & $293N$ \\
    \hline  
\end{tabular}
\label{tab:store1}
\end{table}

\subsection{Flow through porous media}
\label{res:porous}
Next, we analyze several solvers for the simulation of porous media flow. 
We determine optimal parameters for rGCROT(m,k) and GMRES(m) and   
compare the solvers and preconditioners for those parameters. 
Finally, the tuning of the hybrid approach is discussed.

\subsubsection*{Parameter study for rGCROT(m,k) and GMRES(m)}
For a fair comparison, the parameters of rGCROT and 
GMRES(m) are optimized independently for each solver.
For rGCROT, we consider the performance with the Jacobi 
preconditioner while varying the restart frequency ($m$) 
and the dimension of the recycle space ($k$). The parameter $k$ is varied from $10$ to $210$, with increments of $10$, and 
$m$ is varied from $10$ to $80$, with increments of $5$. 
For the other parameters, we take the following values: 
$k-10$ for the number of outer vectors after truncation 
of the outer space; $m/2$ for the number of inner vectors 
to select outer vectors from; $2$ for the number of 
inner vectors selected to extend the outer space; and, $0$
for the number of latest inner vectors kept. 

We obtained the the minimum runtime for rGCROT(10,40), and 
this choice is used throughout this section.

The restart frequency for GMRES(m) with the Jacobi preconditioner
is varied from $15$ to $55$, with increments of $5$, and 
with the SSOR preconditioner from $30$ to $60$, with increments 
of $10$. GMRES(30) with the Jacobi preconditioner and GMRES(50) with the SSOR preconditioner gave the optimal performance in terms of time to solution for the first 10 time steps, and this choice
is used for further comparison.

\subsubsection*{Comparison of solvers and preconditioners}
\label{porous_compare}
\begin{figure}
    \centering
    \includegraphics[width = 12 cm]{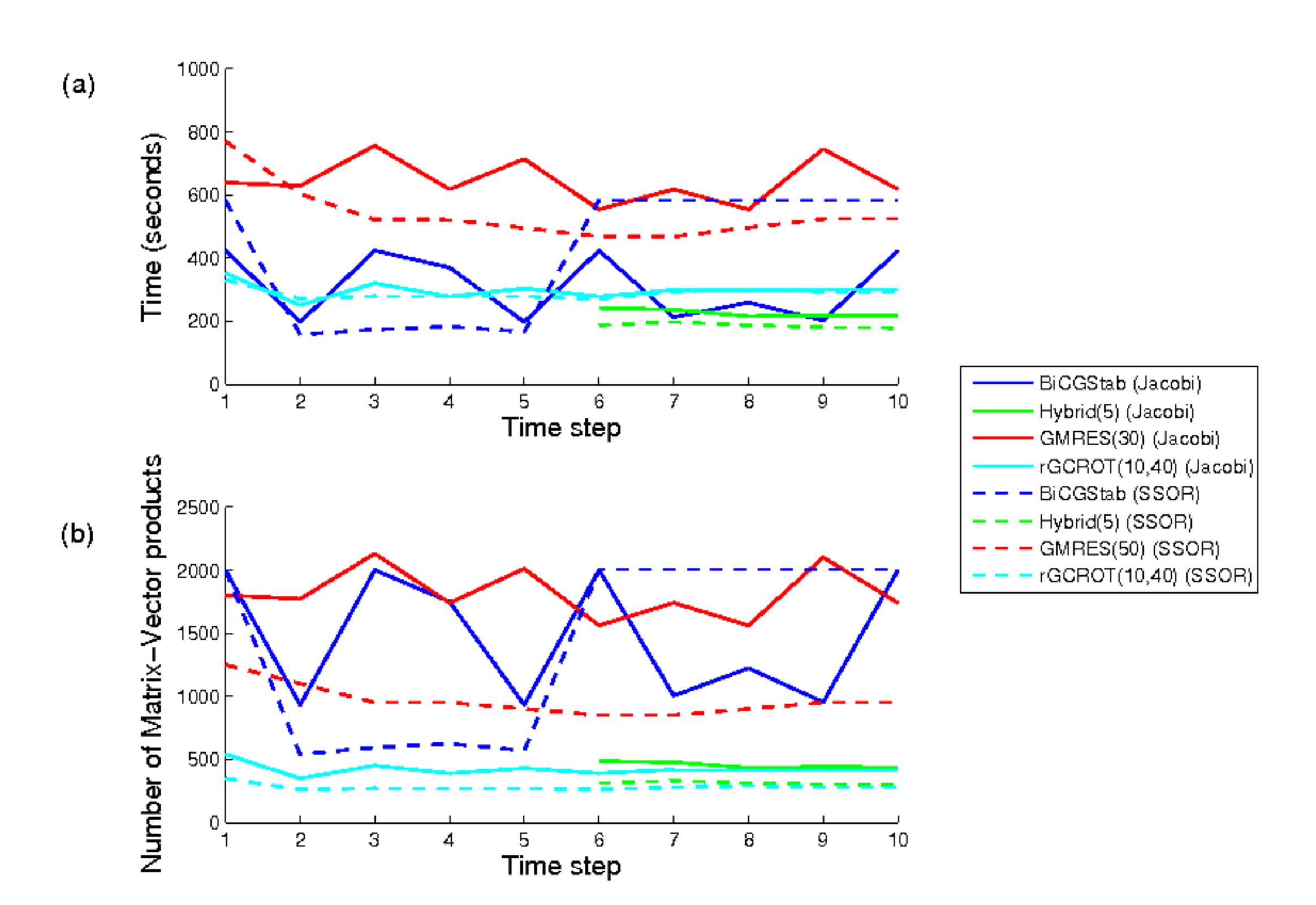}
    \caption{Comparison of solvers using the Jacobi preconditioner (solid lines) and 
    using the SSOR preconditioner (dashed lines): (a) Time to converge (seconds), 
     (b) Number of matrix-vector products to convergence to a 
     relative tolerance of $10^{-10}$.}
    \label{fig:porous1}
\end{figure}

We compare the average number of matrix-vector products and total time taken for the first 10 time steps in Table~\ref{tab:porous1}. 
Note that the average numbers of matrix-vector products for BiCGStab are based on a maximum of 1000 iterations (2000 matrix-vector products) for a time step, and for several time steps this method did not reach the convergence tolerance (see figure~\ref{fig:porous1}). We continued the simulation for the purpose
of comparison.
The average number of matrix-vector products for the hybrid approach
is taken only over the time steps that were 
solved using rBiCGStab. For the hybrid approach, hybrid(n) indicates $n$ time steps with rGCROT(10,40) 
and $(10-n)$ time steps with rBiCGStab.

The results in Table~\ref{tab:porous1} show that the 
hybrid method has the lowest time to solution with both 
preconditioners and requires only a moderate average number 
of matrix-vector products for convergence. 
We also see that rGCROT(10,40) itself is competitive, only losing to
the hybrid solver in run time while having the
lowest average number of matrix-vector products required
for convergence. However, the difference in
matrix-vector products is small, and as rBiCGStab is cheaper in terms 
of orthogonalizations 
it reduces the runtime substantially, especially with the SSOR preconditioner.

Using the SSOR preconditioner reduces the average number of
matrix-vector products for all solvers, indicating that it is the more effective preconditioner for this problem. 
For all solvers except BiCGStab, the total time also decreases with the SSOR preconditioner. 

Figure~\ref{fig:porous1} shows the number of matrix-vector products
to converge and the time to solution per time step for the various solvers.  
This complements Table~\ref{tab:porous1}, which only gives the average results. 
The results for the hybrid solver are given only for the 
time steps solved with rBiCGStab, as the results for the initial steps 
with rGCROT(10,40) are already given. 
The number of iterations remains almost constant for all the 
solvers except for BiCGStab. As can be seen, BiCGStab (with either 
preconditioner) did not converge 
to the required relative tolerance for several time steps. 
Clearly, BiCGStab is not robust enough for this problem. 

Table~\ref{tab:porous2} shows the number of matrix-vector
products for convergence and the solution time for the tenth time step. 
We focus on this time step, because the hybrid approach is tuned based on this time step as discussed in the next subsection. 
The hybrid approach is represented by rBiCGStab, since the tenth time step is solved by rBiCGStab using the recycle space from $5$ initial time steps of rGCROT(10,40).
The rBiCGStab solver with the SSOR preconditioner has by far the lowest time to solution.

\begin{figure} 
    \centering
    \includegraphics[width = 14 cm]{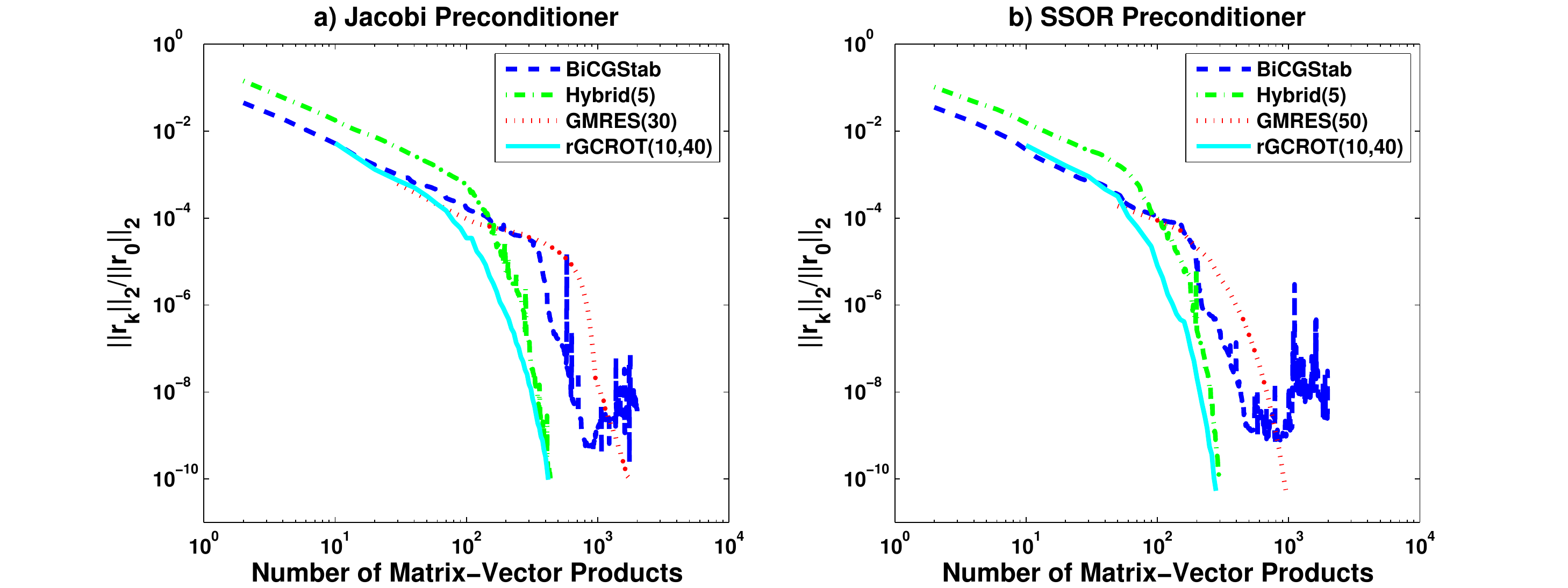}
    \caption{Convergence comparison for the tenth time step using BiCGStab, Hybrid, GMRES and rGCROT solvers with (a) the Jacobi preconditioner and (b) the SSOR preconditioner} 
    \label{fig:porous2}
\end{figure}

Figure~\ref{fig:porous2} plots the convergence of the relative 
residual norm for each solver, $\|r_k\|_2/\|r_0\|_2$,
against the number of matrix-vector products, 
for the tenth time step. The residual norms of both GMRES(m) and rGCROT decrease monotonically, whereas the residual norm of 
BiCGStab shows irregular convergence towards the end
and the method does not converge. The residual norm of
rBiCGStab also displays some (modest) irregular convergence towards the end. 
When comparing the preconditioners, for all solvers the initial residual ($\|r_1\|_2/\|r_0\|_2$) and initial slope are about the same, except for GMRES(m) since the restart frequencies differ.

The $L^2$ norm convergence curves in Figure~\ref{fig:porous2} are given for the 
{\em true} residuals ($b - Ax_i$), except for the rBiCGStab solver, 
where it is given for the {\em updated} residual ($r_{i+1} = s - \omega t$, 
line~15 in Algorithm~4). The update from the outer projection space to 
the solution vector in rBiCGStab is performed only after convergence (see line~17a Algorithm 4). So, the calculation of the 
true residual for every iteration would introduce a
dditional overhead, and is therefore avoided.

The storage costs for the various solvers for the porous media flow 
problem are listed in Table~\ref{tab:store2}. Similar to the discussion in section~\ref{res:turb}, the storage cost for the BiCGStab solver is the lowest,
and it is the highest for the rGCROT algorithm. 
The hybrid approach with rBiCGStab has slightly lower persistent storage requirements than rGCROT(10,40) but higher than for the other solvers. However, it has the
lowest time to solution. So, there is a trade-off between performance and storage,
and the selection of the best solver may depend on the availability of storage  
and data access costs.

\begin{table}[tbp]
    \centering
    \caption{Average number of Matrix-Vector products per time step and total time for 10 time steps of the porous media flow problem for various solvers and preconditioners. The dashes below for GMRES(m) indicate that the convergence for the other m-value was faster (we aim to compare only for the optimal parameters). }
    \begin{tabular}{ l | c | c | c | c }
    \hline  
     & \multicolumn{2}{|c}{Average number of matrix-vector products} & \multicolumn{2}{|c}{Total Time (s)}\\ \cline{2-5}
    Solver & Jacobi & SSOR & Jacobi & SSOR \\
    \hline                       
    BiCGStab & 1480 (max 2000) & 1432 (max 2000) & 3185 & 4202 \\
    Hybrid(5) & 454 & 310 & 2682 & 2076 \\
    GMRES(30) & 1800 & - & 6497 & - \\
    GMRES(50) & - & 950 & - & 5539 \\
    rGCROT(10,40) & 420 & 280 & 3028 & 2943 \\
    \hline  
    \end{tabular}
    \label{tab:porous1}
\end{table}
    
\begin{table}[tbp]
    \centering
    \caption{Convergence comparison for the tenth time step of the porous media flow problem. The dashes below for GMRES(m) indicate that the convergence for the other m-value was faster (we aim to compare only for the optimal parameters).}
    \begin{tabular}{ l | c | c | c | c}
    \hline  
     & \multicolumn{2}{|c}{Number of matrix-vector products} & \multicolumn{2}{|c}{Solution time (s)}\\ \cline{2-5}
    Solver & Jacobi & SSOR & Jacobi & SSOR \\
    \hline                       
    BiCGStab      & (2000) & (2000) & 423.1 & 582.1  \\
    rBiCGStab     & 436    &   298  & 215.9 & 178    \\
    GMRES(30)     & 1740   &    -   & 617.5 & -      \\
    GMRES(50)     & -      &   950  &  -    & 523.2  \\
    rGCROT(10,40) & 420    &   280  & 298.5 & 293.4  \\
    \hline  
    \end{tabular}
    \label{tab:porous2}
\end{table}

\begin{table}[tbp]
\centering
\caption{Solver storage cost comparison for flow through porous media case}
\begin{tabular}{l|c}
    \hline
    Solver & Additional storage cost \\
    \hline  
    BiCGStab & $8N$ \\
    GMRES(50) & $55N$ \\
    rGCROT(10,40) & $93N$\\
    Hybrid(rBiCGStab) & $88N$ \\
    \hline  
\end{tabular}
\label{tab:store2}
\end{table}

\subsubsection*{Hybrid approach tuning}
\label{tuning}
We analyze and tune the performance of the hybrid approach to obtain the 
best parameter choices. This analysis provides the parameter choices 
for the experiments discussed in the previous subsection. 
We vary the dimension of the recycle space, $k$, and the time step, $n$, at which we switch the solver from rGCROT(10,40) to rBiCGStab. We vary $k$ from 
$10$ to $210$, with a step of $10$, while using six time steps to build 
the recycle space using rGCROT. The optimal run time for the hybrid approach is 
also obtained with rGCROT(10,40). For a larger recycle space, the cost of projection becomes too high,  whereas for a  smaller recycle space, the number of
matrix-vector products for convergence becomes too large.

Table~\ref{tab:porous3} lists the time to solution for the tenth time step
and the total time for several $n$ in the hybrid(n) approach. 
We use the Jacobi preconditioner for the tuning. 

If we run 1 or 2 time steps with rGCROT(10,40) and then switch to rBiCGStab, 
the rBiCGStab solver becomes unstable, generating ever larger residual norms. 
This suggests that the recycle space is not yet of good quality.
However, with further initial time steps by rGCROT(10,40), 
the recycle space becomes better, and subsequent systems are
solved faster. 
With the minimum runtime as the criterion, 
the optimal switching point is after 5 time steps of rGCROT(10,40).

\begin{figure}
    \centering
    \includegraphics[width = 10 cm]{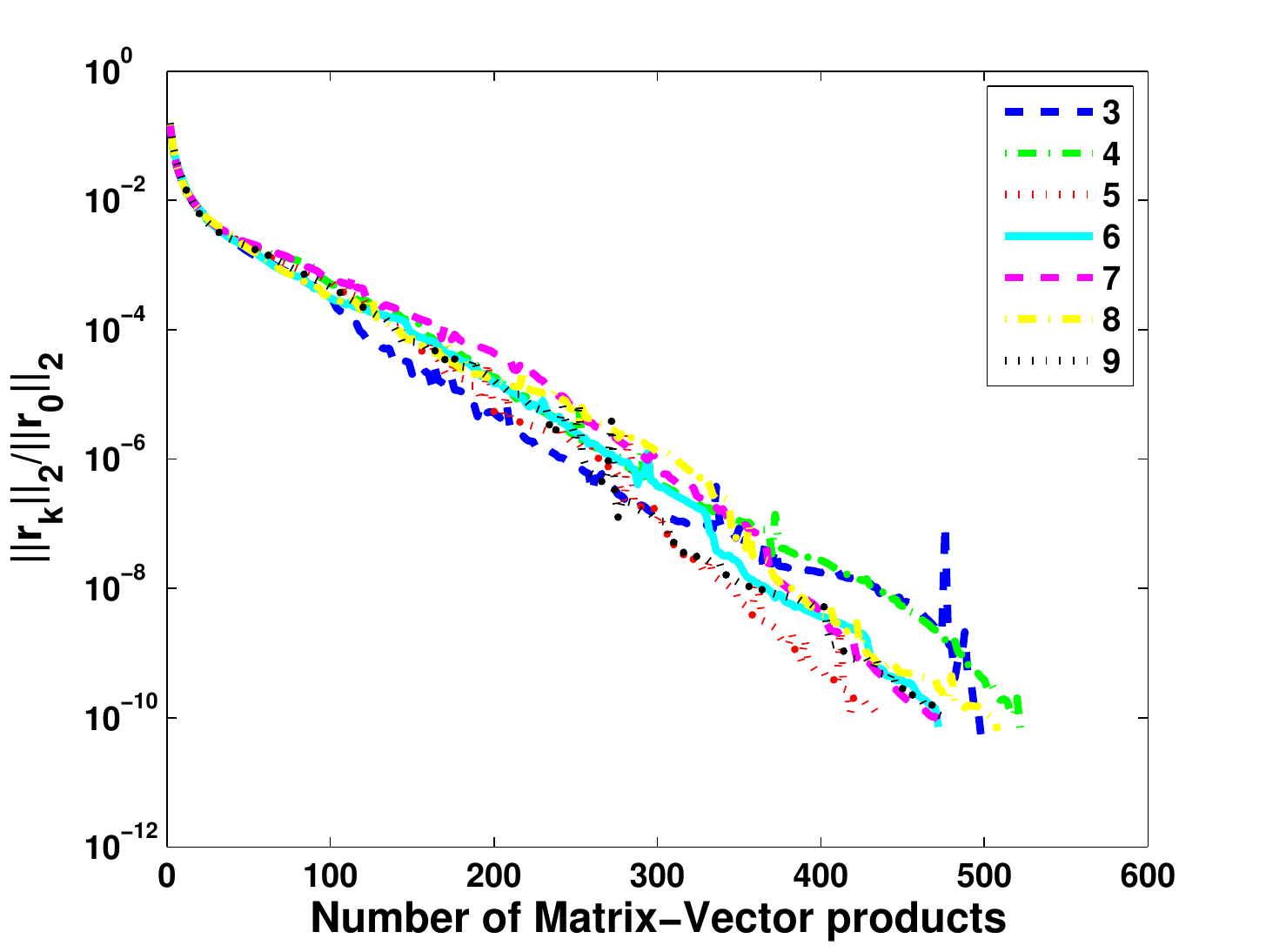}
    \caption{Effect of the number of initial time steps with rGCROT(10,40) on 
      the convergence of the tenth time step with rBiCGStab.}
    \label{fig:porous3}
\end{figure}

The convergence histories for the tenth time step (with rBiCGStab) 
for several n are shown in Figure~\ref{fig:porous3}. The convergence histories for all 
$n > 2$ are fairly close to each other,
suggesting that once a good recurrent subspace is discovered, 
the need to update this space further is minimal. 
These results also show that the choice of time step to switch from 
rGCROT to rBiCGStab is not particularly sensitive, although
care should be taken not to switch too early.
For some problems, it might be useful to update 
the recycle space after a certain number of time steps, 
based on the changes to the system matrix and/or right hand side. 
This was not investigated in the current study.

The hybrid approach is appropriate for problems where the condition number of the preconditioned system matrix is high or the spectrum is otherwise unfavorable, requiring a large number of iterations to converge, 
and there is a sequence of such systems to be solved.
In this case, typically, performance gains from reducing the number of iterations outweighs the overhead per iteration associated with orthogonalization against the recycle space.
Finally, as our results show, the hybrid approach with rBiCGStab 
significantly reduces overall runtime for certain CFD applications 
at the cost of a moderate increase in storage costs.

\begin{table}[tbp]
    \centering
    \caption{Effect of the number of initial time steps with rGCROT(10,40) on 
      the solution time for the hybrid approach.}
    \begin{tabular}{ p{3cm} | p{4cm} | c }
    \hline  
    Number of intial time steps with rGCROT(10,40) &  Solution time for the tenth time step (s) & total time (s) \\
    \hline                       
    1 & unstable & - \\
    2 & unstable & - \\
    3 & 239.5 & 3089 \\
    4 & 248.3 & 3003 \\
    5 & 215.9 & 2682 \\
    6 & 232.8 & 2792 \\
    7 & 233.9 & 2945 \\
    8 & 252.3 & 2928 \\
    9 & 234.4 & 2967 \\
    \hline  
    \end{tabular}
    \label{tab:porous3}
\end{table}    
    
\section{Conclusions and Future Work}
\label{sec5}
The Krylov subspace recycling linear solvers rGCROT and rBiCGStab are tested for two important CFD applications that lead to a sequence of linear systems, and they are compared with the standard solvers GMRES(m) and
BiCGSTab. For a turbulent channel flow problem, the parameter tuned rGCROT solver performs better than both BiCGStab and GMRES(m) in terms of time to solution and number of matrix-vector products to converge. A novel approach to build the recycle space for recycling BiCGStab is proposed and tested for a porous media flow simulation. This hybrid approach has the fastest time to solution 
with a moderate number of matrix-vector products for convergence compared with BiCGSTAB, GMRES(m) and rGCROT.

Several other solver aspects 
provide further important directions of research. 
We should consider a similar study to analyze
which preconditioners provide the best performance, both in terms
of convergence and robustness, including convergence on a
sequence of finer meshes, and in terms of high performance
implementations. In particular, it will be interesting to 
analyze what type of multilevel preconditioners are most 
effective for the CFD problems considered in this paper
and allow high performance on modern architectures,
for example, multiple GPU cards. The optimal interaction between 
Krylov recycling and preconditioning would be an important
aspect of such a study.
Analogously, it would be important to consider how
finer discretizations might affect 
the benefits of recycling.
Of course, such an analysis cannot be separated from 
choices how to precondition for a sequence of finer meshes.

\section{Acknowledgement}
This work was funded in part by the Air Force Office of Scientific Research (AFOSR) Computational Mathematics Program via Grant Number FA9550-12-1-0442 and in part by the National Science Foundation (USA) under Grant Number DMS-1025327 and DMS-1217156. The authors also acknowledge Advanced Research Computing at Virginia Tech for providing computational resources.

Finally, we thank the anonymous reviewers for their careful and
helpful suggestions, which greatly helped us to improve this paper.


\bibliographystyle{elsarticle-num}

\bibliography{references}

\end{document}